# Indefinite quadratic forms and the invariance of the interval in Special Relativity

John H. Elton
Georgia Institute of Technology
elton@math.gatech.edu

**Abstract**.

In this note, a simple theorem on proportionality of indefinite real quadratic forms is proved, and is used to clarify the proof of the invariance of the interval in Special Relativity from Einstein's postulate on the universality of the speed of light; students are often rightfully confused by the incomplete or incorrect proofs given in many texts. The result is illuminated and generalized using Hilbert's Nullstellensatz, allowing one form to be a homogenous polynomial which is not necessarily quadratic. Also a condition for simultaneous diagonalizability of semi-definite real quadratic forms is given.

## 1 Introduction.

In the Special Theory of Relativity, an *event* is a point in space-time whose coordinates with respect to an inertial reference frame correspond to some point $(t, x, y, z)$ in $\mathbb{R}^4$. Coordinates of events in different inertial reference frames are assumed to be connected by linear transformations, based on the assumption of homogeneity and isotropy of space-time. A famous postulate of Einstein is the *universality of the speed of light*: the speed of light in a vacuum is the same in all inertial reference frames, independent of the motion of the source. One can use the postulate of the universality of the speed of light, together with the assumption that changes of coordinates are linear, to determine what changes of coordinates are possible. The idea is to use this postulate to directly show the invariance of a certain quadratic function of the coordinates, which can in turn be used to determine the linear transformations connecting the coordinates (called Lorentz transformations). Defining the Lorentz transformations as the group of linear transformations which leave this quadratic function invariant is geometrically very appealing. To be most satisfying, and not circular, the invariance of the quadratic function should be shown to be a simple and immediate consequence of the postulates; the Lorentz transformations should only then be developed after that.

Suppose points in space-time are specified by $(t, x, y, z)$ in one inertial reference frame $K$, and by $(t', x', y', z')$ in a second inertial reference frame $K'$ whose origin coincides with the first (that is, $t = 0, x = 0, y = 0, z = 0$ in $K$ corresponds to the same event as $t' = 0, x' = 0, y' = 0, z' = 0$ in $K'$). Let a pulse of light be emitted at this common event. Then events on the wave front have coordinates satisfying $x^2 + y^2 + z^2 - c^2 t^2 = 0$ in system $K$, and also $x'^2 + y'^2 + z'^2 - c^2 t'^2 = 0$ in system $K'$, where $c$, the speed of light, is the same in both systems. This is from Einstein's postulate.

In 1966 the author was taking a course in "modern" physics, and remembers being puzzled by the next step taken in the text [7, pg. 58]. The text simply assumed without further ado that $x^2 + y^2 + z^2 - c^2t^2 = x'^2 + y'^2 + z'^2 - c^2t'^2$ for all events (not just those on the wave front of the pulse, when both expressions are zero) and proceeded to use that for a derivation of the form of the Lorentz transformations. Looking in some other texts, we found the same "unconscious" assumption of the *invariance of the interval* $x^2 + y^2 + z^2 - c^2t^2$. In [5, pg. 90 ], it is even stated that "(3.27) $x^2 + y^2 + z^2 - c^2t^2 = 0$ "; "(3.28) $x'^2 + y'^2 + z'^2 - c^2t'^2 = 0$ "; then the amazing statement "...equating lines 3.27 and 3.28, we conclude $x^2 + y^2 + z^2 - c^2t^2 = x'^2 + y'^2 + z'^2 - c^2t'^2$ ". So our confusion remained unresolved for the moment, puzzled by the logic of "things that are equal when zero are always equal" that seemed to be used in these books.

Next semester the author took a course in classical mechanics using the text by J. B. Marion [2]. Appendix G of that book has a demonstration of the invariance of the interval arguing directly from Einstein's postulates, acknowledging the issue that concerned us. (This text is still popular today). Here is the beginning of the proof given in Appendix G, pg. 558, of that book:

(*)    The wave front is described by $x^2 + y^2 + z^2 - c^2t^2 = s^2 = 0$ in $K$ , and $x'^2 + y'^2 + z'^2 - c^2t'^2 = s'^2 = 0$ in $K'$. "...the equations of the transformation that connect the coordinates $(t, x, y, z)$ in $K$ and $(t', x', y', z')$ in $K'$ must themselves be linear. In such a case the quadratic forms $s^2$ and $s'^2$ can be connected by, at most, a proportionality factor: $s'^2 = \kappa s^2$." (It is then shown by further arguments using homogeneity, isotropy and continuity that in fact $\kappa = 1$).

We are of the opinion that the statement above about the reason for the proportionality of the quadratic forms would be misleading to many readers. It is not generally true that if one quadratic form is the result of making a linear change of variables in another quadratic form, and the two quadratic forms have the same zero set, then they must be proportional (even when this zero set has infinitely many points). Here is a somewhat arbitrary example with three variables: Let $s^2 = 2x^2 + 2y^2 + z^2 - 2xz - 2yz$ , and let $s'^2 = 2x'^2 + 2y'^2 + z'^2 - 2x'z' - 2y'z'$ , where $x' = -2x - 2y + z$, $y' = 2y - 2z$, $z' = -2z$ , so the coordinates are connected by a linear transformation. Algebra shows that $s'^2 = 8x^2 + 16y^2 + 10z^2 + 16xy - 16xz - 24yz$ , which is clearly not proportional to $s^2$. Yet both quadratic forms are zero on the same set of points, which is the infinite set $\{z = x + y, x = y\}$, which is apparent after we reveal that actually $s^2 = (x + y - z)^2 + (x - y)^2$, $s'^2 = 10(x + y - z)^2 + 2(x - y)^2 - 4(x + y - z)(x - y)$. For another sort of example (not really related to the statement in Marion but relevant later in this paper), in two variables, let $s^2 = x^2 + y^2 - 2xy$ and $s'^2 = x^2 - y^2$, so $s^2 = 0 \implies s'^2 = 0$, yet these quadratic forms are not even simultaneously diagonalizable. So it would seem the statement about proportionality of the quadratic forms could use further explanation. The author fashioned a proof for himself, but remained puzzled why the books seemed unconcerned about the logical gap.

Fast-forwarding 43 years, we recently had occasion, after not thinking about physics since being an undergraduate, to come upon this topic again. The 1985 text on general relativity by Schutz [6, pg 32] gives a logically correct argument for the proportionality of the quadratic forms in (*). But this does not seem to have been propagated to the community of physics students and textbook writers. From the 2006 relativity text [3], we find on page 10 essentially the same puzzling statements that occurred in the 1964 text [5] mentioned above:

"$c^2t^2 - x^2 - y^2 - z^2 = 0$"; "$c^2t'^2 - x'^2 - y'^2 + z'^2 = 0$"; "These are equal, so $c^2t^2 - x^2 - y^2 - z^2 = c^2t'^2 - x'^2 - y'^2 + z'^2$".

And we have evidence, from the Physics Forum, that indeed other physics students are still finding themselves confused by exactly the same thing that we found unexplained so long ago! The answers we saw given by other students in the Forum were unfortunately not correct and were essentially on the level of the "unconscious" proofs of some of those texts, along with some rather arrogant statements about the students who didn't understand the "proofs" they saw in their books. See

www.**physicsforums**.com/archive/index.php/t-115451.html

which archives this amusing discussion.

So we decided *this time* to fill in the gap, for the benefit of others who might be confused, by stating and proving a more general but very simple result about indefinite quadratic functions that settles the matter. This result about containment of zero sets suggests a more general result, proved using Hilbert's Nullstellensatz. Also we prove a simple result about simultaneous diagonalization of semidefinite quadratic forms and containment of zero sets.

## 2        A theorem on indefinite quadratic forms.

A function $q : \mathbb{R}^n \to \mathbb{R}$ is a *real quadratic form* if there is a symmetric bilinear function $\tilde{q} : \mathbb{R}^n \times \mathbb{R}^n \to \mathbb{R}$ such that $q(\mathbf{x}) = \tilde{q}(\mathbf{x}, \mathbf{x})$. In matrix language, this means there is a symmetric $n \times n$ matrix $Q = \left[ Q_{ij} \right]$ of real numbers such that

$q(x_1, \dots x_n) = q(\mathbf{x}) = \sum_{i=1}^{n} \sum_{j=1}^{n} Q_{ij} x_i x_j$, i.e., $q(\mathbf{x}) = \mathbf{x}^t Q \mathbf{x}$ for $\mathbf{x} \in \mathbb{R}^n$. The elements of the

matrix $Q$ are the components of $\tilde{q}$ in the standard basis.

A real quadratic form $q$ is *indefinite* if it takes both positive and negative values; this is equivalent to the matrix $Q$ having at least one positive eigenvalue and at least one negative eigenvalue. See [4] for example, or any book on linear algebra.

For a real quadratic form $q$, define $Z_q = \{\mathbf{x} \in \mathbb{R}^n : q(\mathbf{x}) = 0\}$; this is the *zero set* of $q$.

**Theorem 1.** Let $q$ be an indefinite real quadratic form on $\mathbb{R}^n$, and let $r$ be a real quadratic form on $\mathbb{R}^n$ such that $Z_q \subset Z_r$; that is, $q(\mathbf{x}) = 0 \implies r(\mathbf{x}) = 0$. Then $r$ is proportional to $q$; that is, there exists a real number $\alpha$ such that $r(\mathbf{x}) = \alpha q(\mathbf{x})$ for all $\mathbf{x}$. If $\alpha$ is not zero, then $r$ is also indefinite and has the same zero set as $q$.

*Proof*: There exists a basis $\{\mathbf{v}_1, ... \mathbf{v}_n\}$ for $\mathbb{R}^n$ such that the matrix $Q_{ij} = \tilde{q}(\mathbf{v}_i, \mathbf{v}_j)$ representing $q$ in this basis is diagonal, with only 1's, -1's and 0's on the diagonal, and $Q_{ii} = 1$ *for* $1 \le i \le k$; $Q_{ii} = -1$ *for* $k+1 \le i \le k+m$; $Q_{ii} = 0$ *for* $k+m+1 \le i \le n$; and $Q_{ij} = 0$ *for* $i \ne j$. The numbers k and m here are unique: $k$ is the number of positive eigenvalues and $m$ is the number of negative eigenvalues of any matrix representing $q$ (Sylvester's law of inertia; see [4, pg. 202]). Since $q$ is indefinite, $k > 0$ and $m > 0$. So without loss of generality, in the proof which follows we will just assume that $Q$ is a diagonal matrix with k ones and m negative-ones and the rest (if any) zeroes on the diagonal, in order, as described above. (In the application to invariance of the interval which motivated this discussion, Q is already of this form, but we wanted to treat the general case). Let $R$ be the symmetric matrix representing $r$ in this basis.

The idea is to make judicious choices of points where $q$ is zero, to conclude that $R$ must also be diagonal, and that the on-diagonal elements of $R$ are a common multiple of those of $Q$.

To that end, let j be an integer such that $k+1 \le j \le k+m$. Let $\mathbf{x}$ have components $x_1 = 1, x_j = 1$, and all other components zero. Then $q(\mathbf{x}) = Q_{11}x_1^2 + Q_{jj}x_j^2 = 1 - 1 = 0$, so $r(\mathbf{x}) = R_{11}x_1^2 + R_{jj}x_j^2 + 2R_{1j}x_1x_j = R_{11} + R_{jj} + 2R_{1j} = 0$, by hypothesis. Now change the sign of the j$^{th}$ component of $\mathbf{x}$ so that $x_j = -1$ but leave the other components of $\mathbf{x}$ unchanged; then $q(\mathbf{x}) = 0$ still, so $r(\mathbf{x}) = R_{11}x_1^2 + R_{jj}x_j^2 + 2R_{1j}x_1x_j = R_{11} + R_{jj} - 2R_{1j} = 0$ also. These two equations together imply $R_{1j} = 0$ and then $R_{jj} = -R_{11}$. Then for $1 < i \le k$, using i in place of 1 in the argument above shows $R_{ii} = -R_{jj} = R_{11}$, and $R_{ij} = 0$.

Next let j be an integer (if any) such that $m+k+1 \le j \le n$. First let $\mathbf{x}$ be the vector with $x_j = 1$ and all other components zero. Then q($\mathbf{x}$) = $Q_{jj}$ = 0, so r($\mathbf{x}$) = $R_{jj}$ = 0. Next let $1 < i \le k$, $k+1 \le l \le k+m$, and let $\mathbf{x}$ be the vector with components $x_i = 1, x_l = 1, x_j = 1$, and all other components zero. Then $q(\mathbf{x}) = Q_{ii} + Q_{ll} + Q_{jj} = 1 - 1 + 0 = 0$, so $r(\mathbf{x}) = R_{ii} + R_{ll} + R_{jj} + 2R_{il} + 2R_{ij} + 2R_{lj} = 2R_{ij} + 2R_{lj} = 0$ also. Changing $\mathbf{x}$ so that $x_i = -1$ and otherwise unchanged leads to $-2R_{ij} + 2R_{lj} = 0$. This implies that $R_{ij} = 0$, and then $R_{lj} = 0$.

Suppose that $k \ge 2$. Let $1 \le i < j \le k$, $k+1 \le l \le k+m$, and let $\mathbf{x}$ be the vector with components $x_i = 3, x_j = 4, x_l = 5$, and all other components zero. Then $q(\mathbf{x}) = Q_{ii}x_i^2 + Q_{jj}x_j^2 + Q_{ll}x_l^2 = 9 + 16 - 25 = 0$, so

$$r(\mathbf{x}) = R_{ii}x_i^2 + R_{jj}x_j^2 + R_{ll}x_l^2 + 2R_{ij}x_ix_j + 2R_{il}x_ix_l + 2R_{jl}x_jx_l = R_{11}(9+16-25) + 2R_{ij}(12) = 0$$

also (note we have already shown that $R_{il} = R_{jl} = 0$ and the proportionality of the diagonal elements). This proves $R_{ij} = 0$. Similarly, if $m \geq 2$, the corresponding off-diagonal terms of $R$ are zero.

This completes the proof that $R = R_{11}Q$, and the proof of the theorem. $\square$

## 3    An alternate proof using Hilbert's Nullstellensatz, and a stronger result.

The containment of zero sets in the hypothesis of theorem 1 suggests Hilbert's Nullstellensatz [1, pg. 254], of importance in algebraic geometry. We can also prove Theorem 1 using this theorem rather than using diagonalization and bases as we did above; and although the proof above is certainly simple enough, there is some insight to be gained by this alternate proof, and a more general result can be proved this way as well. The Nullstellensatz concerns zero sets of ideals in the ring of polynomials in several variables over an algebraically closed field. For our application the ideal in question will be simply the principal ideal generated by a single polynomial q. If the reader is not familiar with ideal theory and the Nullstellensatz, it will not matter because we shall use only the following immediate consequence of Hilbert's theorem:

If $q(\mathbf{x})$ and $r(\mathbf{x})$ are complex polynomials in n variables such that $\mathbf{x} \in \mathbb{C}^n$ and $q(\mathbf{x}) = 0 \Rightarrow r(\mathbf{x}) = 0$, then $r^p(\mathbf{x}) = q(\mathbf{x})s(\mathbf{x})$ for some polynomial $s(\mathbf{x})$ and positive integer $p$. If $q$ is square-free (that is, the irreducible factors of $q$ occur only to the first power), $p$ can be taken to be one.

In theorem 1, $q$ and $r$ are quadratic forms with real coefficients, $q$ is indefinite, and the *real* zeroes of $q$ are assumed to be zeros of $r$ by hypothesis. Unless the rank of $q$ is 2 (it can't be less than 2 if $q$ is indefinite), q is easily seen to be irreducible (the zero set of $q$ isn't just the union of two subspaces when the rank is three or more, so $q$ couldn't be the product of two linear factors); and even if $q$ were reducible, it would be the product of linear factors to the first power and therefore square-free, so we can always take $p$ to be one in our application. Thus all we need to do is to show that, as a consequence of the indefiniteness of $q$, the complex zeroes of $q$ are also zeroes of $r$, and the conclusion of theorem 1 will follow from the Nullstellensatz, since the degrees of $q$ and $r$ being two requires $s$ to be constant.

To that end, suppose $q(\mathbf{x}+i\mathbf{y}) = 0$ for some $\mathbf{x}, \mathbf{y} \in \mathbb{R}^n$, so $q(\mathbf{x}) - q(\mathbf{y}) = 0$ and $\tilde{q}(\mathbf{x}, \mathbf{y}) = 0$. If $q(\mathbf{x}) = 0$ (hence $q(\mathbf{y}) = 0$) then $q(\mathbf{x}+\mathbf{y}) = 0$, so $r(\mathbf{x}) = r(\mathbf{y}) = r(\mathbf{x}+\mathbf{y}) = 0$, which implies $\tilde{r}(\mathbf{x}, \mathbf{y}) = 0$ and so $r(\mathbf{x}+i\mathbf{y}) = 0$.

Suppose then that $q(\mathbf{x}) > 0$ (the opposite case would be handled similarly); by rescaling assume $q(\mathbf{x}) = 1$. Since $q$ is indefinite, there is $\mathbf{u} \in \mathbb{R}^n$ such that $q(\mathbf{u}) < 0$. Let $\mathbf{w} = \mathbf{u} - \tilde{q}(\mathbf{u}, \mathbf{x})\mathbf{x} - \tilde{q}(\mathbf{u}, \mathbf{y})\mathbf{y}$, so $\tilde{q}(\mathbf{w}, \mathbf{x}) = 0$ and $\tilde{q}(\mathbf{w}, \mathbf{y}) = 0$ (a "Gram-Schmidt" construction). Now $q(\mathbf{w}) = \tilde{q}(\mathbf{w}, \mathbf{u}) = q(\mathbf{u}) - \tilde{q}(\mathbf{u}, \mathbf{x})^2 - \tilde{q}(\mathbf{u}, \mathbf{y})^2 < 0$. By rescaling we may assume that $q(\mathbf{w}) = -1$, $\tilde{q}(\mathbf{w}, \mathbf{x}) = 0$ and $\tilde{q}(\mathbf{w}, \mathbf{y}) = 0$. Thus $q(\mathbf{w} + \alpha\mathbf{x} + \beta\mathbf{y}) = -1 + \alpha^2 + \beta^2 = 0$ whenever $\alpha^2 + \beta^2 = 1$, so by hypothesis $r(\mathbf{w} + \alpha\mathbf{x} + \beta\mathbf{y}) = r(\mathbf{w}) + \alpha^2 r(\mathbf{x}) + \beta^2 r(\mathbf{y}) + 2\alpha\tilde{r}(\mathbf{w},\mathbf{x}) + 2\beta\tilde{r}(\mathbf{w},\mathbf{y}) + 2\alpha\beta\tilde{r}(\mathbf{x},\mathbf{y}) = 0$ also for $\alpha^2 + \beta^2 = 1$. Taking $\alpha = \pm1, \beta = 0$, we conclude that $\tilde{r}(\mathbf{w},\mathbf{x}) = 0$, and similarly $\tilde{r}(\mathbf{w},\mathbf{y}) = 0$. Then choosing $\alpha = 2^{-1/2}, \beta = \pm\alpha$, we conclude that $\tilde{r}(\mathbf{x},\mathbf{y}) = 0$. Then choosing $\alpha = 1, \beta = 0$ and $\alpha = 0, \beta = 1$, we see that $r(\mathbf{x}) = r(\mathbf{y})$, and thus $r(\mathbf{x} + i\mathbf{y}) = 0$, concluding the proof.

This proof from the Nullstellensatz is perhaps slightly cleaner than the first proof of theorem 1. But also one can prove more this way, with a little more work. The quadratic form $r$ is a polynomial in n variables in which each term has degree 2. In general, a polynomial in n variables for which each term has the same degree d is called a *homogeneous polynomial* of degree d.

**Theorem 1a.** Suppose $r$ is a homogeneous real polynomial in n variables, not necessarily a quadratic form, with the other hypotheses of theorem 1 unchanged. Then $q$ is a factor of $r$; that is, $r(\mathbf{x}) = q(\mathbf{x})s(\mathbf{x})$ for some polynomial $s(\mathbf{x})$.

*Proof*: We only need to show that any complex zeroes of $q$ are zeroes of $r$. Suppose $q(\mathbf{x} + i\mathbf{y}) = 0$, and suppose that $q(\mathbf{x}) > 0$. As above, we can assume that $q(\mathbf{x}) = q(\mathbf{y}) = 1$, $\tilde{q}(\mathbf{x},\mathbf{y}) = 0$, and there is $\mathbf{w}$ such that $q(\mathbf{w}) = -1$, $\mathbf{w}$ is $q$-orthogonal to $\mathbf{x}$ and $\mathbf{y}$, and $q(\mathbf{w} + \alpha\mathbf{x} + \beta\mathbf{y}) = 0$, so $r(\mathbf{w} + \alpha\mathbf{x} + \beta\mathbf{y}) = 0$ also, whenever $\alpha^2 + \beta^2 = 1$. Suppose that $r$ has even degree $2m$. Now $r(\mathbf{w} + \alpha\mathbf{x} + \beta\mathbf{y})$ is a polynomial of degree 2m in the variables $\alpha, \beta$ which is zero on the unit circle $\alpha^2 + \beta^2 = 1$. We may write $r(\mathbf{w} + \alpha\mathbf{x} + \beta\mathbf{y}) = \sum_{j+k \le 2m} \alpha^j \beta^k c(j,k)$ where the indices $j$ and $k$ are nonnegative.

$c(j,k)$ is the sum of terms involving the coefficients of the polynomial r and the components of $\mathbf{w}, \mathbf{x}$ and $\mathbf{y}$ in general, but for terms of maximum degree 2m in the variables $\alpha, \beta$, $\mathbf{w}$ is not involved. By changing the signs of $\alpha$ and $\beta$ separately, and then together, we see that for $\alpha^2 + \beta^2 = 1$, $\sum_{j+k \le 2m, j \text{ odd}, k \text{ even}} \alpha^j \beta^k c(j,k) = 0$,

$\sum_{j+k \le 2m, j \text{ even}, k \text{ odd}} \alpha^j \beta^k c(j,k) = 0$, $\sum_{j+k \le 2m, j \text{ and } k \text{ odd}} \alpha^j \beta^k c(j,k) = 0$, and $\sum_{j+k \le 2m, j \text{ and } k \text{ even}} \alpha^j \beta^k c(j,k) = 0$.

Consider the last expression above (with both indices even), which can be rewritten with a change of indices as $\sum_{j+k \le m} (\alpha^2)^j (1-\alpha^2)^k c(2j, 2k) = 0$ for $\alpha^2 \le 1$. This is a polynomial of degree 2m in $\alpha$; the coefficient of the highest power term $\alpha^{2m}$ must be zero because of the constancy of the polynomial on an infinite set, so $\sum_{j \le m} (-1)^{m-j} c(2j, 2m-2j) = 0$. Next consider the next-to-last expression (with both indices odd) which can be rewritten $\alpha\beta \sum_{j+k \le m-1} (\alpha^2)^j (1-\alpha^2)^k c(2j+1, 2k+1) = 0$, so for $0 < \alpha^2 < 1$, $\sum_{j+k \le m-1} (\alpha^2)^j (1-\alpha^2)^k c(2j+1, 2k+1) = 0$  Setting the coefficient of the highest power term in the polynomial in $\alpha$ (which occurs when $k = m-j-1$) to zero, we get

$\sum_{j \le m-1} (-1)^{m-j-1} c(2j+1, 2m-2j-1) = 0$.

Since $r$ is homogeneous of degree $2m$,

$r(\mathbf{x} + i\mathbf{y}) = \sum_{j \le m} (-1)^{m-j} c(2j, 2m-2j) + i \sum_{j \le m-1} (-1)^{m-j-1} c(2j+1, 2m-2j-1)$, because

$i^{2m-2j} = (-1)^{m-j}$ and $i^{2m-2j-1} = (-1)^{m-j-1} i$. The results just proved show this is zero, completing the proof of the theorem when the degree of $r$ is even and $q(\mathbf{x}) > 0$.

Now suppose r has odd degree $2m-1$. Then

$r(\mathbf{w} + \alpha\mathbf{x} + \beta\mathbf{y}) = \sum_{j+k \le 2m-1} \alpha^j \beta^k c(j, k)$, and this breaks into four sums equaling zero on the unit circle as before, depending on the parities of the indices. Consider first consider the sum corresponding to $j$ odd and $k$ even; this can be rewritten

$\alpha \sum_{j+k \le m-1} (\alpha^2)^j (1-\alpha^2)^k c(2j+1, 2k) = 0$ for $\alpha^2 \le 1$. Setting the coefficient of the highest power to zero gives $\sum_{j \le m-1} (-1)^{m-j-1} c(2j+1, 2m-2j-2) = 0$. Now consider the sum corresponding to $j$ even and $k$ odd, which can be rewritten

$\beta \sum_{j+k \le m-1} (\alpha^2)^j (1-\alpha^2)^k c(2j, 2k+1) = 0$, so $\sum_{j+k \le m-1} (\alpha^2)^j (1-\alpha^2)^k c(2j, 2k+1) = 0$ for

$\alpha^2 < 1$ which implies $\sum_{j \le m-1} (-1)^{m-j-1} c(2j, 2m-2j-1) = 0$. But

$r(\mathbf{x} + i\mathbf{y}) = \sum_{j \le m-1} (-1)^{m-j-1} c(2j+1, 2m-2j-2) + i \sum_{j \le m-1} (-1)^{m-j-1} c(2j, 2m-2j-1)$, so this is

zero, and the proof is concluded for the case $r$ is of odd degree and $q(\mathbf{x}) > 0$.

The case when $q(\mathbf{x}) < 0$ is handled in a similar way.

Finally, if $q(\mathbf{x} + i\mathbf{y}) = 0$ and $q(\mathbf{x}) = 0$, then since $q(\mathbf{y}) = 0$ and $\tilde{q}(\mathbf{x}, \mathbf{y}) = 0$,

$q(\mathbf{x} + \alpha\mathbf{y}) = 0$ and thus $r(\mathbf{x} + \alpha\mathbf{y}) = 0$ for all real numbers $\alpha$. Now $r(\mathbf{x} + \alpha\mathbf{y})$ is a polynomial in $\alpha$ which is identically zero, so all its coefficients are zero, and this clearly implies that $r(\mathbf{x} + i\mathbf{y}) = 0$, which concludes the proof of theorem 1a. $\square$

.

## 4    Simultaneous diagonalization of quadratic functions.

Theorem 1 implies a result on simultaneous diagonalizability: if $q$ is an indefinite real quadratic form on $\mathbb{R}^n$ and $r$ is a real quadratic form on $\mathbb{R}^n$ such that $Z_q \subset Z_r$, then $q$ and $r$ are simultaneously diagonalizable (meaning there is a basis in which the matrices representing $q$ and $r$ are both diagonal).

However, if $q$ is a semi-definite real quadratic form on $\mathbb{R}^n$ (*semi-definite* means $q(\mathbf{x}) \ge 0$ for all $\mathbf{x}$ or $q(\mathbf{x}) \le 0$ for all $\mathbf{x}$), and $r$ is a real quadratic form on $\mathbb{R}^n$ such that $Z_q \subset Z_r$, then $q$ and $r$ are not necessarily simultaneously diagonalizable. For an example in $\mathbb{R}^2$, let $q(\mathbf{x}) = (x - y)^2$ and let $r(\mathbf{x}) = x^2 - y^2$; this example (already mentioned in the introduction) satisfies the conditions and it is easy to see these are not simultaneously diagonalizable.

But if $q$ and $r$ are both assumed semi-definite, there is a similar (and similarly easy) result on containment of zero sets implying simultaneous diagonalizability.

**Theorem 2.** Let $q$ and $r$ be semi-definite real quadratic forms on $\mathbb{R}^n$ such that $Z_q \subset Z_r$. Then $r$ and $q$ are simultaneously diagonalizable.

*Proof*: Without loss of generality assume they are both positive semi-definite. First observe that the zero sets are subspaces: Let $q(\mathbf{x}) = 0$ and $q(\mathbf{y}) = 0$; then $q(a\mathbf{x} + b\mathbf{y}) = a^2 q(\mathbf{x}) + b^2 q(\mathbf{y}) + 2ab\tilde{q}(\mathbf{x}, \mathbf{y}) = 2ab\tilde{q}(\mathbf{x}, \mathbf{y}) \geq 0$ for all real numbers $a, b$ implies $\tilde{q}(\mathbf{x}, \mathbf{y}) = 0$, so $q(a\mathbf{x} + b\mathbf{y}) = 0$. This is quite different from the indefinite case where the zero sets are cones and not subspaces.

There is a subspace $M$ such that $\mathbb{R}^n = M \oplus Z_q$ (choose any basis for $Z_q$ and extend it to a basis for $\mathbb{R}^n$, and $M$ is the span of those added-on basis vectors). Let $\mathbf{x} = \mathbf{y} + \mathbf{z}$ with $\mathbf{y} \in M$ and $\mathbf{z} \in Z_q$. Then $q(\mathbf{y} + \alpha\mathbf{z}) = q(\mathbf{y}) + 2\alpha\tilde{q}(\mathbf{y}, \mathbf{z}) + \alpha^2 q(\mathbf{z}) = q(\mathbf{y}) + 2\alpha\tilde{q}(\mathbf{y}, \mathbf{z}) \geq 0$ for all real $\alpha$ implies $\tilde{q}(\mathbf{x}, \mathbf{y}) = 0$, so $q(\mathbf{y} + \mathbf{z}) = q(\mathbf{y})$. Similarly, $r(\mathbf{y} + \mathbf{z}) = r(\mathbf{y})$ for $\mathbf{y} \in M$ and $\mathbf{z} \in Z_q$, since $Z_q \subset Z_r$.

Thus $q$ and $r$ may be considered as positive semi-definite quadratic forms on $M$, and in fact $q$ is positive definite on $M$, because if $\mathbf{y} \in M$ and $q(\mathbf{y}) = 0$, then $\mathbf{y} \in Z_q$ by definition, so $\mathbf{y} = 0$. By a well-known theorem[4, pg. 218], this implies $q$ and $r$ are simultaneously diagonalizable on $M$, and they are then simultaneously diagonalizable on $\mathbb{R}^n = M \oplus Z_q$, with zeroes on the diagonal corresponding to the basis vectors for $Z_q$. $\square$

## 5    Application to the proof of invariance of the interval.

Suppose the coordinates $\mathbf{x} = (t, x, y, z)$ in $K$ and $\mathbf{x}' = (t', x', y', z')$ in $K'$ are connected by a linear transformation, so $\mathbf{x}' = L\mathbf{x}$ for some $4 \times 4$ matrix $L$. Let $q(\mathbf{x}) = -c^2 t^2 + x^2 + y^2 + z^2 = \mathbf{x}^t Q\mathbf{x}$, where $Q$ is the diagonal matrix with diagonal entries $(-c^2, 1, 1, 1)$. Let $r(\mathbf{x}) = -c^2 t'^2 + x'^2 + y'^2 + z'^2 = (L\mathbf{x})^t QL\mathbf{x} = \mathbf{x}^t (L^t QL)\mathbf{x}$, so $r(\mathbf{x}) = \mathbf{x}^t R\mathbf{x}$, where $R = L^t QL$. Now $q$ is indefinite, and $r(\mathbf{x}) = 0$ precisely when $q(\mathbf{x}) = 0$, from (*) above. So the conditions of Theorem 1 are in force, so we may conclude that $r$ is proportional to $q$, which is equivalent to the statement from (*) that we wanted to prove, namely that $s'^2$ is proportional to $s^2$.

**Acknowledgment**. We would like to thank Michael Loss for suggesting looking at Hilbert's Nullstellensatz for a connection with the topic in this note.